\documentclass[final,onefignum,onetabnum]{siuro210301}

\usepackage{lipsum}
\usepackage{amsfonts}
\usepackage{graphicx}
\usepackage{epstopdf}
\usepackage{algorithmic}
\ifpdf
  \DeclareGraphicsExtensions{.eps,.pdf,.png,.jpg}
\else
  \DeclareGraphicsExtensions{.eps}
\fi

\usepackage{enumitem}
\setlist[enumerate]{leftmargin=.5in}
\setlist[itemize]{leftmargin=.5in}

\newsiamremark{remark}{Remark}
\newsiamremark{hypothesis}{Hypothesis}
\crefname{hypothesis}{Hypothesis}{Hypotheses}
\newsiamthm{claim}{Claim}

\headers{Applications of AD in Image Registration}{W. Watson, C. Cherry, and R. Lang}

\title{Applications of Automatic Differentiation in Image Registration}

\author{Warin Watson\thanks{Department of Mathematics and Statistics, Colorado Mesa University
  (\email{wdwatson2@mavs.coloradomesa.edu}).}
\and Cash Cherry\thanks{Department of Applied Mathematics and Statistics, Colorado School of Mines 
  (\email{ccherry@mines.edu}).}
\and Rachelle Lang\thanks{Department of Mathematics, University of Wisconsin-Madison (\email{rklang@wisc.edu}).}}

\dedication{\small\textit{Project advisor: Lars Ruthotto\thanks{Department of Mathematics, Emory University}}}

\usepackage{amsopn}

\newcommand{\bfb}{{\bf b}}

\newcommand{\bfx}{{\bf x}}

\newcommand{\bfd}{{\bf d}}

\newcommand{\bfr}{{\bf r}}
\newcommand{\bff}{{\bf f}}
\newcommand{\bfv}{{\bf v}}
\newcommand{\bfw}{{\bf w}}

\newcommand{\bfom}{{\boldsymbol{\omega}}}

\newcommand{\CT}{{\cal T}}
\newcommand{\CR}{{\cal R}}

\DeclareMathOperator*{\argmin}{arg\,min}

\newcommand{\RR}{\mathbb{R}}
\newcommand{\vy}{\vec{y}}
\newcommand{\vx}{\vec{x}}

\newcommand{\sfA}{\mathsf{A}}

\newcommand{\sfD}{\mathsf{D}}

\newcommand{\sfI}{\mathsf{I}}

\newcommand{\sfK}{\mathsf{K}}
\newcommand{\sfL}{\mathsf{L}}

\newcommand{\sfO}{\mathsf{O}}

\ifpdf
\hypersetup{
  pdftitle={Applications of Automatic Differentiation in Image Registration},
  pdfauthor={W. Watson, C. Cherry, and R. Lang}
}
\fi

\begin{document}

\maketitle

\begin{abstract}
We demonstrate that automatic differentiation (AD), which has become commonly available in machine learning frameworks, is an efficient way to explore ideas that lead to algorithmic improvement in multi-scale affine image registration and affine super-resolution problems. In our first experiment on multi-scale registration, we implement an ODE predictor-corrector method involving a derivative with respect to the scale parameter and the Hessian of an image registration objective function, both of which would be difficult to compute without AD. Our findings indicate that exact Hessians are necessary for the method to provide any benefits over a traditional multi-scale method; a Gauss-Newton Hessian approximation fails to provide such benefits. In our second experiment, we implement a variable projected Gauss-Newton method for super-resolution and use AD to differentiate through the iteratively computed projection, a method previously unaddressed in the literature. We show that Jacobians obtained without differentiating through the projection are poor approximations to the true Jacobians of the variable projected forward map and explore the performance of other approximations in the problem of super-resolution. By addressing these problems, this work contributes to the application of AD in image registration and sets a precedent for further use of machine learning tools in this field.
\end{abstract}

\section{Introduction}

In the context of medical imaging, image registration is the problem of finding a reasonable transformation $\vy$ to align a template image $\CT$ with a reference image $\CR.$ In the optimization framework, we find a transformation by minimizing the sum of a distance function $\mathcal{D}[\CT \circ \vy, \CR]$ and a regularization term \cite{Modersitzki2009}.  Approaches to find optimal transformations in a parameterized set of admissible transformations can be found either by numerical optimization techniques, or by solving a nonlinear PDE derived from the optimality condition \cite{Fischer_2008}. Our strategy, following \cite{Modersitzki2009}, is to parameterize the transformations by a vector and numerically solve for the optimal transformation as a finite-dimensional unconstrained optimization problem. 

The problem of super-resolution, which we properly introduce in \S 4, seeks to reconstruct an unknown high resolution reference image from a sequence of unregistered low resolution templates \cite{chung2006numerical}. Solution techniques for this problem implement on algorithms for image registration while simultaneously reconstructing the high resolution reference image. We use the same strategy as in registration to solve the super-resolution problem by using a finite dimensional approximation to the problem. However, due to the simultaneous reconstruction of the reference image and the registration of the templates, the objective has a separable structure amenable to a variable projection \cite{varpro} approach. 

Leveraging the PyTorch \cite{PyTorch} library, we consider two novel applications of automatic differentiation (AD) to this framework. We apply a predictor-corrector method to perform multi-scale image registration and improve existing variable projection methodology to solve a super-resolution problem. The prediction step of the predictor-corrector method requires a derivative with respect to a scaling parameter in the interpolation, and variable projection involves differentiation through an inexact iterative least squares solve. Much of the theory for these problems exists in the literature, but perhaps due to the difficulty of some required derivatives, they do not exist in standard image registration packages like FAIR \cite{Modersitzki2009}. The code and data required generate the figures used in this text and to reproduce results can be found at our GitHub repository \cite{repo}.

\section{Image Registration Background}
We provide background on the optimization setting used in our work, and the capabilities of Automatic Differentiation.

\subsection{Numerical Optimization Framework}

Assume our images are supported on a box-shaped subset $\Omega \subset \RR^d$ for $d=2$ or $3.$ We consider an enumerated set $\bfom = \{ \vx_1, \dots, \vx_n \} \subset \Omega$ of points with even spacing $h_{x}$ and $h_{y}$ in each dimension, on which we have evaluations of the reference and template images. (We do \emph{not} have evaluations for points not on the grid, which we will have to address.) Using this data, we would like to be able to compute the distance $\mathcal{D}[\CT \circ \vy, \CR]$ of our choice, which in this work will be the least-squares distance metric
\begin{equation}\label{LS_cont}
    \mathcal{D}[\CT \circ \vy, \CR] = \int_\Omega ( \CT(\vy(\vx)) - \CR(\vx) )^2 dV.
\end{equation}
As an approximation to the integral \eqref{LS_cont}, we use the quadrature rule $\| \CT( \vy(\bfom) ) - \CR(\bfom) \|_2^2 h_x h_y,$ where $\CT( \vy(\bfom) )$ and $ \CR(\bfom)$ are taken to be vectors of evaluations on $\bfom$ of $\CT \circ \vy$ and $\CR,$ respectively. To compute the evaluations $\CT(\vy(\bfom))$, we must evaluate the template image at points $\vy(\vx_i)$ which will usually not be contained in the grid $\bfom.$ Further, as we want to use first and second order optimization techniques to minimize a function containing this term, we need a twice continuously differentiable interpolant. To accomplish this, we use the spline interpolation scheme as described in \cite{Modersitzki2009}. (This hypothetically allows for \emph{exact} quadrature of the interpolants; however, it would be more difficult to code and it's not clear it would make a difference in the registration quality.)

The final thing to discretize is the transformation $\vy.$ The most expressive ``non-parametric'' approach would be to solve for every value in $\vy(\bfom)$ separately. However, to avoid such a high dimensional optimization problem, we use a parameterization in many fewer parameters $\bfw \in \RR^p$, henceforth writing $\vy(\vec{\cdot}; \bfw)$. A number of standard parameterizations are covered in \cite{Modersitzki2009}. In this work, we consider only affine transformations for simplicity. Our goal is to evaluate and demonstrate the behavior of the proposed methods, rather than to solve difficult image registration problems.

Our discretized image registration loss function has the general form
\begin{equation}\label{eq:discrete_loss}
    J(\bfw) = \| \CT(\vy(\bfom; \bfw)) - \CR(\bfom) \|_2^2 h_x h_y + \lambda^2 S ( \bfw ),
\end{equation}
where $S ( \bfw )$ is a regularization term, and $\lambda > 0$ is the corresponding regularization parameter. Note that cubic splines are twice continuously differentiable, so the least squares distance $\| \CT(\vy(\bfom; \bfw)) - \CR(\bfom) \|_2^2 h_x h_y$ is twice continuously differentiable with respect to $\bfw$ as long as $\vec{y}(\bfom; \bfw)$ is. Note also that the third derivative of a cubic spline with fixed data exists and is bounded \emph{almost everywhere} (as cubic splines are piecewise cubic polynomials on intervals); so as long as $\vec{y}(\bfom; \bfw)$ also has bounded third derivatives with respect to $\bfw,$ the Hessian of the least squares loss is also Lipschitz continuous. If we furthermore choose regularization term satisfying this condition, then \eqref{eq:discrete_loss} is twice continuously differentiable function on $\RR^p$ with a Lipschitz continuous Hessian, thus we can rigorously expect Newton's method, Gauss-Newton, and gradient descent to converge to local minimizers with their respective rates of convergence \cite{NoceWrig06}.

\subsection{Automatic Differentiation}

In the past, a major difficulty in implementing first and second order numerical optimization schemes on complicated objective functions has been the analytic computation of derivatives. However, with automatic differentiation (AD), the implementation of such methods are simplified by no longer having to hand code analytical derivatives, which is known to be time-consuming and error-prone \cite{margossian2019review}. 

To demonstrate the practical efficiency of automatic differentiation (AD), Table~\ref{tab:computation_times} reports average evaluation times for computing various derivatives under two transformation models: an affine transformation with 6 parameters, and a neural ODE-based transformation with 50 parameters~\cite{chen2018neural, SUN2024103249, wu2022nodeo}. 

We define the \emph{forward function} as the map \(\bfw \mapsto \CT(\vy(\bfom; \bfw))\), which evaluates the template image at coordinates determined by \(\bfw\), prior to loss computation. The full loss function \(J(\bfw)\) follows \eqref{eq:discrete_loss}, incorporating squared differences and quadrature. All derivatives are computed using Functorch~\cite{functorch2021}, which provides composable AD transforms similar to those in JAX~\cite{schoenholz2020jax}.

While the theoretical cost of computing gradients and Hessians scales as $\mathcal{O}(n)$ and $\mathcal{O}(n^2)$, respectively, modern AD frameworks such as Functorch exploit graph-level optimizations, vectorization, and batched computation to reduce overhead~\cite{functorch2021}. We observe that while the gradient and Jacobian costs scale moderately, the Hessian cost increases substantially for the 50-parameter case. This supports the practical distinction in cost between Newton and quasi-Newton methods for higher-dimensional transformations.

\begin{table}[H]
\centering
\caption{
Empirical cost of computing derivatives using AD. “Factor” is the ratio of runtime relative to the corresponding loss or forward function evaluation. For example, computing the gradient took 2.5× the time of a loss evaluation in the affine case. Each measurement was averaged over 100 runs. Reported values are mean runtime (\(\mu\)) ± standard deviation (\(\sigma\)) in milliseconds.
}

\label{tab:computation_times}
\vspace{0.5em}
\begin{tabular}{l|rr|rr}
 & \multicolumn{2}{c|}{6 parameters} & \multicolumn{2}{c}{50 parameters} \\
\cline{2-3} \cline{4-5}
 & Time (ms) ($\mu \pm \sigma$) &  Factor & Time (ms) ($\mu \pm \sigma$) & Factor \\
 \hline
Loss Function     & $2.171 \pm 0.259$ & $1\times$ & $5.012 \pm 2.168$ & $1\times$ \\
Gradient          & $5.422 \pm 0.409$ &  $2.5\times$  & $13.65 \pm 1.018$ &  $2.7\times$  \\
Hessian           & $12.44 \pm 1.798$ &  $5.7\times$  & $120.0 \pm 8.737$ &  $24\times$  \\
\hline
Forward Function  & $1.054 \pm 0.076$ & $1\times$ & $3.204 \pm 0.149$ & $1\times$ \\
Jacobian          & $11.59 \pm 1.518$ &  $11\times$  & $54.31 \pm 3.240$ &  $17\times$  \\
\end{tabular}
\end{table}

Many (real-valued) linear operators associated with these problems (particularly for the Super Resolution problem) are cheap to apply to a given vector, but expensive to store explicitly. For the purpose of solving normal equations associated with them, it is necessary to compute products with the adjoints of these operators. This is done easily using an explicitly stored matrix form by taking the transpose, but is not straightforward when we cannot store the matrix. 

Forward-mode automatic differentiation computes the derivative $\frac{\partial f (\bfx + h\bfv)}{\partial h}$ with respect to any one dimensional perturbation $\bfv$ of the input space at the cost of a single forward evaluation. For a matrix-vector product $f ( \bfx ) = \sfA \bfx$, such a derivative gives the product $\bfv^\top \sfA,$ for which $\sfA^\top \bfv = ( \bfv^\top \sfA )^\top,$ the desired adjoint-vector product. Crucially, the matrix $\sfA$ need not be explicitly constructed, as long as the function for computing its action is written in an automatically differentiable way. 
 
\section{Multi-scale Methods}

\begin{figure}[]
    \centering
    \includegraphics[width=.95\linewidth]{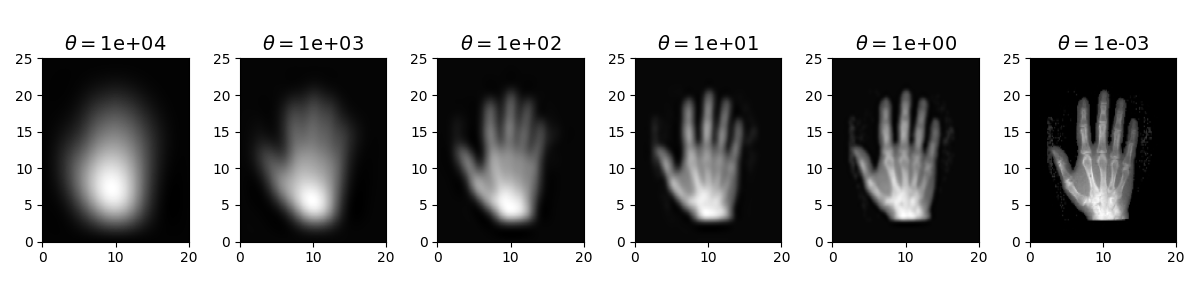}
    \caption{$\theta = 0$ corresponds to a regularly interpolated image with no blurring, while increasing $\theta$ increases the amount of blurring.}
    \label{fig:scaling}
\end{figure}

One of the main difficulties of solving image registration problems is that the objective functions \eqref{eq:discrete_loss} resulting from real images are often non-convex and contain many local minima due to fine-scale details. Thus, by solving a very down-sampled or blurred registration problem, the minimum obtained is closer to the global minimizer of the desired problem \cite{Modersitzki2009}. This leads to an iterative approach, where the solutions to easier versions of the registration problem are iteratively used as initializations for harder versions. There are two such commonly used approaches in registration: The \emph{multi-level} approach (which we don't consider in this work) down-samples the images, and the \emph{multi-scale} approach blurs the images. Figure \ref{fig:scaling} illustrates the blurring effect that the \emph{scaling parameter}, $\theta$, has on the interpolation of an image. For further explanation and examples of the multi-level and multi-scale techniques, the reader is directed to \cite{Modersitzki2009}. 

An ODE, derived in Section \ref{sec:homo_deriv}, describes how the optimal set of transformation parameters changes with respect to $\theta$. One could follow the path of a minimizer from a large value of $\theta$ down to a small value of $\theta$ purely by numerically solving the ODE, but in our experience, this approach is computationally expensive. Instead, we opt for a predictor-corrector method, where the ODE is used to refine the initial guess at each scale, and local minimization is employed to correct errors at that scale. We distinguish these multi-scale methods from traditional multi-scale methods due to the utilization of the ODE. Additionally, we distinguish these methods from homotopy methods \cite{osti_876373} since the scaling is applied to the interpolation of the images rather than directly to the objective function.

\subsection{Derivation of Multi-scale ODE}\label{sec:homo_deriv}

Consider a minimizer $\bfw^*$ for a fixed $\theta$ image registration problem, defined by:
    \begin{equation*}
        \bfw^*(\theta) \in \argmin_{\bfw \in \RR^p} J ( \bfw; \theta ).
    \end{equation*}
As $J$ is twice continuously differentiable, and $\bfw^*$ is a local minimizer, it must satisfy the first-order optimality condition
\begin{equation}
    \label{eq:condition}
        \nabla_{\bfw} J( \bfw^*(\theta); \theta ) = 0 .
\end{equation}
Differentiating \eqref{eq:condition} with respect to $\theta$,
$$
\frac{d}{d\theta}\nabla_{{\bfw}} J( \bfw^*(\theta); \theta ) = \nabla^2_{{\bfw}} J( \bfw^*(\theta); \theta ) \frac{d{\bfw}^*}{d\theta} + \frac{d}{d\theta} \nabla_{\bfw} J ( \bfw^*(\theta); \theta ) = 0.
$$
It follows that
\begin{equation}
\label{eq:homo_ode}
    \frac{d}{d\theta}{\bfw}^*(\theta) = -\bigg( \nabla^2_{{\bfw}} J (\bfw^*(\theta); \theta)\bigg)^{-1} \left(\frac{d}{d\theta} \nabla_{\bfw} J (\bfw^*(\theta); \theta)\right).
\end{equation}
Each step of this ODE solve requires an inverse Hessian-vector product, so it is as feasible to implement as Newton's method for optimization.

\subsection{Predictor-Corrector Method}

The approach involves taking larger steps in $\theta$ when solving \eqref{eq:homo_ode} (i.e., using coarser discretizations for the ODE). The error introduced by these larger steps is then corrected through local minimization. In practice, Newton's method for the minimization provides rapid local convergence \cite{NoceWrig06}, and allows us the ability to re-use computed Hessians when taking a step in the ODE solve. 

To take a step from a coarse scaling parameter $\theta_n$ to a finer scaling parameter $\theta_{n+1} < \theta_n$ starting from a minimizer $\bfw^*(\theta_{n+1})$, a step
\begin{equation}
\label{eq:predictor}
    \bfw^*_{\rm pred} = \bfw^*(\theta_{n}) + (\theta_{n} - \theta_{n+1})\left(\frac{d}{d\theta}\bfw^*(\theta_{n})\right)
\end{equation}
of forward Euler's method is taken to predict the location of $\bfw^*(\theta_{n+1}),$ and is then corrected by a local minimization procedure to obtain the true $\bfw^*(\theta_{n+1}).$ 

Even with the predictor-corrector method, negative curvature in the local minimization problem still poses a problem. The forward Euler predictor step \eqref{eq:predictor} explicitly uses the Hessian inverse, which can lead to computational instability if the Hessian is indefinite. Pure Newton steps also struggle in this scenario, since directly using indefinite Hessians can produce non-descent directions. To address this, we opt for using SciPy’s implementation of an exact trust-region Newton solver \cite{conn2000trust} for the local minimization problem, which avoids explicit Hessian inversion by solving a trust-region subproblem through eigen-decomposition~\cite{2020SciPy-NMeth}. A line-search Newton method is another viable alternative~\cite{NoceWrig06}. After convergence of the local minimization, we reuse the Hessian computed at the accepted iterate in the subsequent forward Euler predictor step, avoiding additional Hessian evaluations.

Note that the predictor-corrector method reduces to a traditional multi-scale method when the forward Euler step of the ODE \eqref{eq:predictor} is not taken. This suggests that to prefer the predictor-corrector method over the traditional multi-scale method, you would need to ensure the cost of computing \eqref{eq:predictor} does not outweigh the cost of performing minimization to get the same $\bfw^*$. This would vary based on the problem and behavior of the optimization algorithm and its implementation. 

\subsection{Results}

Two results are shown in this section. First, the predictor-corrector method is shown to work for a non-regularized problem where registering from a single scale does not. This motivates the method as having a similar use case as the traditional multi-scale method. Then, for the predictor-corrector method, the prediction is shown to fail when using approximated Hessians via Gauss Newton, but succeeds when using Hessians computed using AD, which are exact with respect to the discretized objective~\eqref{eq:discrete_loss}. An implementation of Levenberg-Marquardt is used for Gauss Newton, while SciPy's Exact Trust Region \cite{2020SciPy-NMeth} is used for Newton.

\subsubsection{Single Scale vs Predictor-Corrector}

\begin{figure}[]
    \centering
    \includegraphics[scale=0.7]{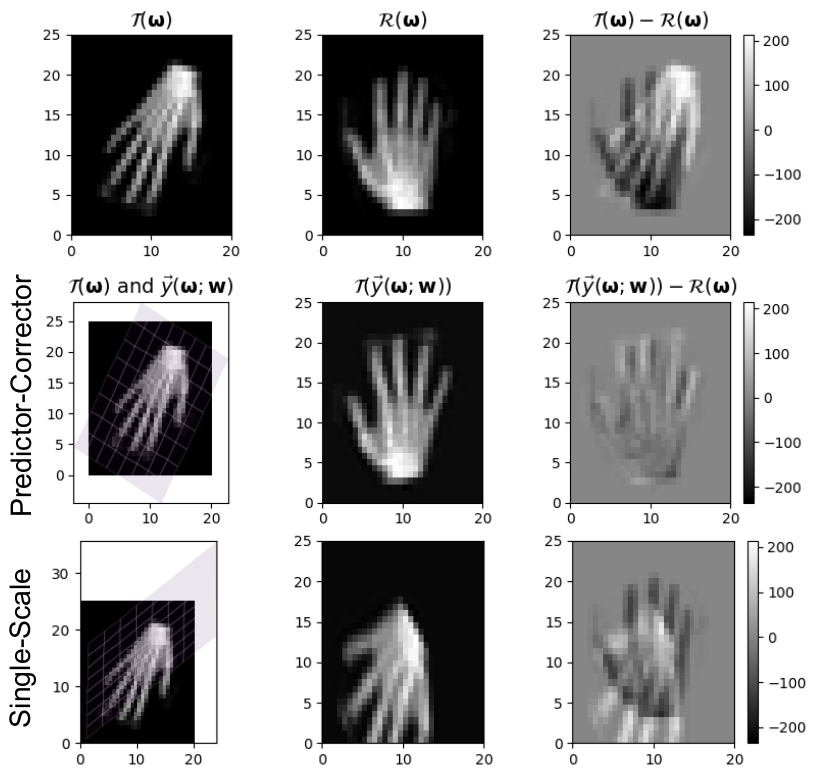}
    \caption{Template and Reference are shown in the first row. The predictor-corrector method is used to register the images in the second row. Attempting to register at the finest scale alone results in a bad local minimizer as shown in the third row.}
    \label{fig:singlescale}
\end{figure}

Consider an affine, non-regularized image registration problem with template and reference shown in the first row of Figure \ref{fig:singlescale}. As shown in the figure, some image registration problems benefit from using the predictor-corrector method, as opposed to optimizing from a single scale, in order to avoid getting stuck in a local minimum. 

\subsubsection{Approximated Hessian vs Exact Hessian}

Gauss-Newton is often implemented over a Newton method to optimize non-linear least squares problems as a way to avoid expensive and/or complicated Hessian computations. Let us define the residual function as $r(\bfw) = \CT(\vy(\bfom; \bfw)) - \CR(\bfom),$ and consider the un-regularized least squares loss function
$$
J ( \bfw ) = \| r(\bfw)  \|_2^2 .
$$
The exact Hessian of such a function is given by 
$$
\nabla^2 J ( \bfw ) = \nabla r ( \bfw )^\top \nabla r ( \bfw ) + \sum_{i=1}^n r_i (x) \nabla^2 r_i (x),
$$
where $\nabla r ( \bfw )$ is the Jacobian of $r ( \bfw )$. Utilizing the exact Hessian in Newton's method yields a quadratic rate of convergence, but can be expensive. At the cost of relaxing the guarantee to only sub-linear convergence, Gauss Newton avoids the computational expense of computing second derivatives by approximating the Hessian using solely the term $\nabla r ( \bfw )^\top \nabla r ( \bfw )$. However, the method can approach quadratic convergence and often does in practice \cite{NoceWrig06}. With AD, both methods are simple to implement, though the computational expense of computing the Hessian for Newton's method is still greater.

To evaluate the effectiveness of the prediction at each scale, i.e., the step of forward Euler \eqref{eq:predictor}
, we compare the loss and gradient at the predicted point to the loss and gradient without the prediction. This comparison effectively compares the predictor-corrector method with the traditional multi-scale method. We define the relative loss difference as $(J(\bfw^*) - J(\bfw_{\rm pred}^*) )/J(\bfw^*)$, and the relative gradient norm difference as $(\| \nabla J(\bfw^*) \| - \| \nabla J(\bfw_{\rm pred}^*) \| ) /$ $ \| \nabla J(\bfw^*) \|$ to give measures of how helpful the predictions are. 

\begin{figure}[]
    \centering
    \includegraphics[width=0.8\linewidth]{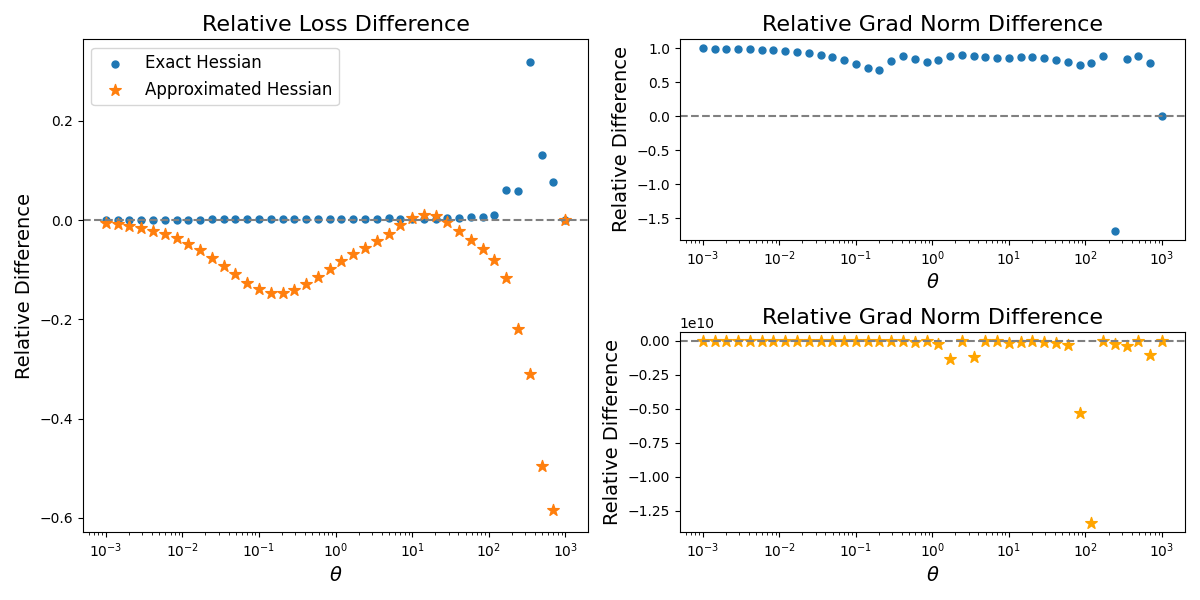}
    \caption{Relative loss and gradient norm, both plotted against $\theta$ on a logarithmic scale. Exact Hessians are used for blue, approximated Hessians are used for orange. All relative gradient norm difference values using approximated Hessians are less than $-31000$. }
    \label{fig:diffplots}
\end{figure}

This experiment is illustrated in Figure \ref{fig:diffplots}. A positive value indicates a successful reduction in loss or gradient, whereas a negative value indicates an increase. The results show that using exact Hessians via Newton's method tends to decrease the loss and gradient, while using approximated Hessians via Gauss-Newton tends to increase the loss and gradient. Therefore, using exact Hessians is necessary for the predictor-corrector method to provide any benefits over a traditional multi-scale method. 

\section{Coupled Methods for Super-Resolution}
We follow the framework developed in \cite{chung2006numerical}. Interpolation of the reference is \emph{linear} in the reference intensity vector $\bff^{(0)}$, thus, we can mathematically model our transformations $\vy^{(j)}$ using matrices $\sfI^{(j)} (\bfw).$ The notation $\sfI$ with no argument is used for the identity matrix, which is consistent with viewing the identity as an application of a trivial transformation. Note that the interpolation matrix is \emph{non-linear} in its dependence on the transformation parameter vector $\bfw.$ Our high resolution transformed templates $\bff^{(j)} \in \RR^n, j = 1, \dots, q$ are thus given by the transformations $\sfI^{(j)} (\bfw) \bff^{(0)}$ of an initial reference image $\bff^{(0)}$. Our observed template images $\bfd^{(j)} \in \RR^m, j = 0, 1 \dots, q$ are modeled as $\bfd^{(j)} = \sfK\bff^{(j)},$ where $\sfK \in \RR^{m \times n}$ is a linear operator which lowers the image resolution ($m < n$), in our case by pooling, that is, averaging the pixels in $k \times k$ regions to reduce resolution by $k$ times in each dimension. Note that our modeling assumes that the template $\bfd^{(0)}$ is registered with the unknown reference $\bff^{(0)}$, as the data can only ever constrain $\bff^{(0)}$ up to a transformation; registering $\bfd^{(0)}$ with $\bff^{(0)}$ eliminates this non-uniqueness and saves some computational overhead. We further let $h^{\bff}_x, h^{\bff}_y, h^{\bfd}_x$ and $h^{\bfd}_y$ be the corresponding grid sizes for the template and reference grids.

In the problem of super-resolution, we aim to reconstruct a high-resolution reference image, $\mathbf{f}^{(0)}$ from a vector of low resolution template images $\mathbf{d} = [{\mathbf{d}^{(0)}}^\top, \dots, {\mathbf{d}^{(q)}}^\top]^\top$, which are assumed to be unregistered, down-sampled versions of $\mathbf{f}^{(0)}$. Our aim in super-resolution is to register multiple templates $\mathbf{d}^{(j)}$ to an unknown reference $\mathbf{f}^{(0)},$ while simultaneously reconstructing $\mathbf{f}^{(0)}$ from the unregistered templates.

\subsection{Super-Resolution Framework}

Representing the template intensities in a single vector $\bfd,$ our full model for the super resolution data is given by
$$
\bfd = (\sfI \otimes \sfK) \sfI (\bfw) \bff^{(0)}, \quad \quad \sfI (\bfw) = \begin{bmatrix}
    \sfI \\ \sfI^{(1)} (\bfw) \\ \vdots \\ \sfI^{(q)} (\bfw)
\end{bmatrix}, \quad \quad \bfd = 
\begin{bmatrix}
     \bfd^{(0)} \\ \vdots \\ \bfd^{(q)}
\end{bmatrix}.
$$
We could approach this as a regularized least squares problem in the two variables $\bff^{(0)}$ and $\bfw$: 
$$
J ( \bfw, \bff^{(0)} ) = \| (\sfI \otimes \sfK) \sfI (\bfw) \bff^{(0)} - \bfd \|^2 h^{\bfd}_x h^{\bfd}_y +  \lambda_{\bff}^2 S_{\bff^{(0)}} (\bff^{(0)}).
$$
As before, we don't regularize the transformation parameters $\bfw.$ Our regularization term on the reference image $\bff^{(0)}$ is a finite difference approximation to the semi-norm $\int_{\Omega} \| \nabla \CR^{(0)} (\vx) \|^2 dV,$ where $\CR^{(0)}$ denotes the continuous reference image approximated by the vector $\bff^{(0)}.$ This is given by
$$S_{\bff^{(0)}} (\bff^{(0)}) = \| \sfL \bff^{(0)} \|^2 h^{\bff}_x h^{\bff}_y, \quad \quad \sfL = \begin{bmatrix}
    h^{\bff}_x \sfD \otimes \sfI\\
    \sfI \otimes h^{\bff}_y \sfD
\end{bmatrix},\quad \quad \sfD = 
\begin{bmatrix}
    -1 & 1 &\\
    & \ddots & \ddots &\\
    && -1 & 1
\end{bmatrix}. $$

Our joint super resolution objective function is written
\begin{equation}\label{eq:superres_joint}
    J ( \bfw, \bff^{(0)} ) = \| (\sfI \otimes \sfK) \sfI (\bfw) \bff^{(0)} - \bfd \|^2 + \lambda_{\bff}^2 \| \sfL \mathbf{f}^{(0)} \|^2.
\end{equation}
 
It should be noted that this quadratic regularizer is not ideal, as it penalizes high wavenumber content in the reconstruction that may be part of the true reference image. This shows up in our examples as ``checkerboard'' patterns in the residual. To move beyond this, nonlinear regularization methods such as total variation can be used with a ``Linearize-And-Project'' approach \cite{Herring2018}, however, like in \cite{chung2006numerical}, we stick to a quadratic regularization term as it simplifies both the theory and practice of variable projection.

\subsection{Variable Projection Theory}
We begin by noticing that \eqref{eq:superres_joint} is in the separable form addressed by the variable projection method \cite{varpro,golub2003separable,español2024convergenceanalysisvariableprojection}. To make this even more clear, we write it in the stacked form
\begin{equation*}
    J ( \bfw, \bff^{(0)} ) = \left \| \begin{bmatrix}
    \sqrt{h^{\bfd}_x h^{\bfd}_y} \ (\sfI \otimes \sfK) \sfI (\bfw)  \\ \lambda_{\bff}\sqrt{h^{\bff}_x h^{\bff}_y} \ \sfL 
    \end{bmatrix} \bff - \begin{bmatrix}
         \sqrt{h^{\bfd}_x h^{\bfd}_y} \ \bfd \\ \sfO
    \end{bmatrix} \right \|^2 = \| \sfA_{\lambda_{ \bff}} (\bfw) \bff^{(0)} - \bfb \|^2.
\end{equation*}
In this case, the optimal $\bff^{(0)}$ as a function of $\bfw$ is given by $\bff^{(0)} ( \bfw ) = \sfA_{\lambda_{\bff}} (\bfw)^{\dagger} \bfb.$ Thus, our variable projected objective function is given by
\begin{equation}\label{eq:superres_varpro}
    \tilde{J} (\bfw) = \| \sfA_{\lambda_{\bff}} (\bfw) \bff^{(0)} (\bfw) - \bfb \|^2 
\end{equation}
A natural approach to minimize \eqref{eq:superres_varpro} is Gauss-Newton, in which case we need to compute the Jacobian of the term $\sfA_{\lambda_{\bff}} (\bfw) \bff^{(0)} (\bfw).$ Applying the product rule,

\begin{equation}
   \frac{\partial}{\partial \bfw} \left( \sfA_{\lambda_{\bff}} (\bfw) \bff^{(0)} (\bfw) \right) = \frac{\partial \sfA_{\lambda_{\bff}} (\bfw)}{\partial \bfw} \bff^{(0)} (\bfw) + \sfA_{\lambda_{\bff}} (\bfw) \frac{\partial \bff^{(0)} (\bfw)}{\partial \bfw}.
   \label{eq:Jacobianofresidual}
\end{equation}

Computing the Jacobian $\frac{\partial \bff^{(0)} (\bfw)}{\partial \bfw}$ of the projection $\bff^{(0)} (\bfw)$ requires differentiation through the least squares solution.  An exact formula for such a derivative based on implicit differentiation through an exact solution exists \cite{varpro}, but it is expensive to compute. Previous work in variable-projected Gauss-Newton has employed approximations to avoid this cost—either through low-rank SVD approximations \cite{nnvarpro} or by entirely neglecting the Jacobian term (effectively setting $\frac{\partial \bff^{(0)} (\bfw)}{\partial \bfw} = 0$) as done in \cite{chung2006numerical}.

In an automatic differentiation framework, we can ignore the projection term by turning gradient tracking off during the least squares solve. However, by choosing to track gradients through the least squares solution, we can include this term at the cost of increased memory overhead. It has been shown that automatic differentiation through the Conjugate Gradient (CG) algorithm gives correct gradients in the case of an inexact solve \cite{CGdiff1}, and in the case of an exact solve \cite{CGdiff}, which we will use in our solutions. A recent preprint \cite{otherdiff} explores automatic differentiation through other least-squares algorithms.

\subsection{Implementation and Results}

\begin{figure}[]
    \centering
    \includegraphics[width=0.95\linewidth]{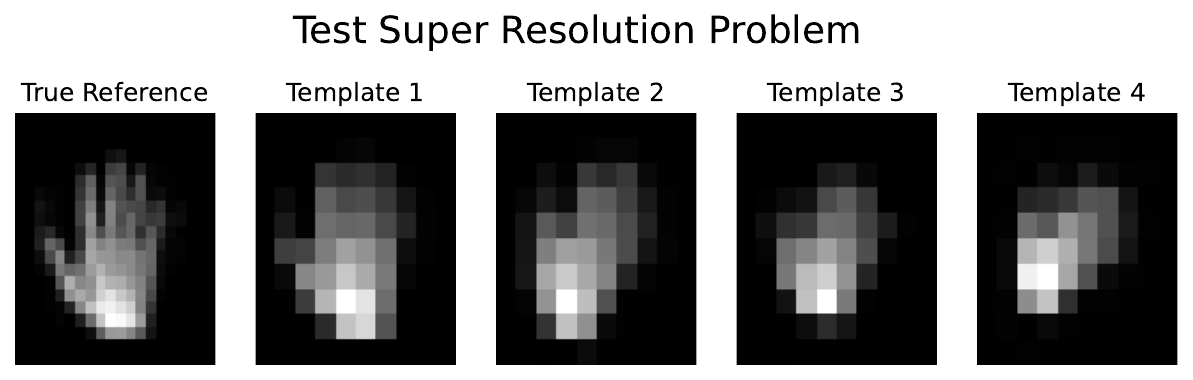}
    \caption{Our test problem for super resolution. The reference image is $20 \times 20$ pixels, and the $4$ templates are $10 \times 10$ pixels. The problem is well-determined, as we are to infer $400$ parameters from $400$ pixels.}
    \label{fig:inexact_supperes_test}
\end{figure}

\begin{figure}[]
    \centering
    \includegraphics[width=0.7\linewidth]{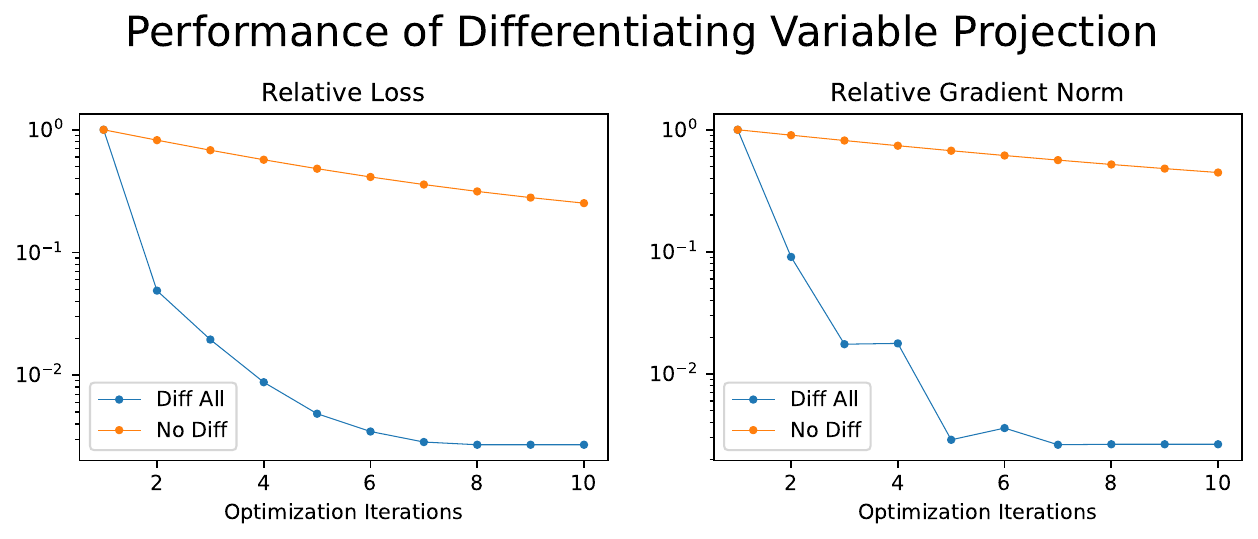}
    \caption{The performance of 10 iterations of variable-projected Gauss-Newton on the problem in Figure \ref{fig:inexact_supperes_test}. Differentiating through all CG iterations provides a faster reduction in relative loss and relative gradient norm. }
    \label{fig:original-performance-plot}
\end{figure}

\begin{figure}[]
    \centering
    \includegraphics[width=0.6\linewidth]{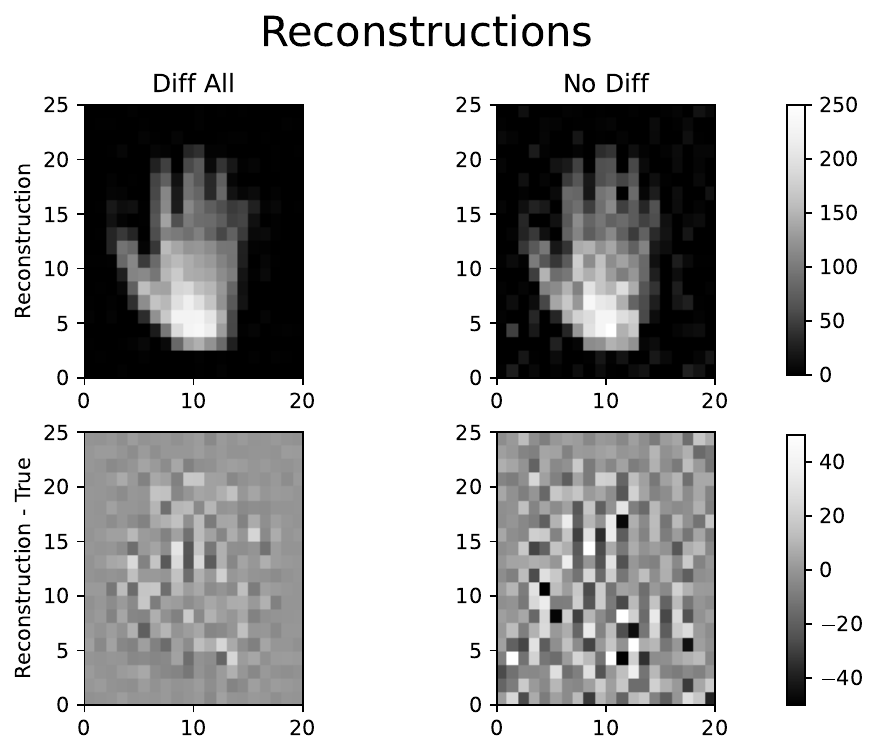}
    \caption{The final reconstructions from the performance test comparing no differentiation against full differentiation through every iteration of CG. }
    \label{fig:original-reconstruction-plot}
\end{figure}

We consider the small super-resolution problem shown in Figure \ref{fig:inexact_supperes_test} and note a significant benefit in reconstruction quality and convergence rate of the Gauss-Newton method when using AD to differentiate through all iterations of CG (see Figure \ref{fig:original-performance-plot} and Figure \ref{fig:original-reconstruction-plot}). There are $18$ registration parameters to estimate in this problem, and we found that computing the Jacobian of the variable-projected objective function (Equation \eqref{eq:superres_varpro}) using forward-mode automatic differentiation takes approximately six to ten times longer than evaluating the objective function itself, which is in line with our expectations. The final relative reconstruction error when differentiating all iterations of CG is $10.7\%$, while ignoring the Jacobian of the projection results in a higher relative reconstruction error of $25.61\%$.

At each optimization step, the Jacobian \eqref{eq:Jacobianofresidual} is computed using forward mode AD. Its important to note that in forward mode, memory consumption scales with the size of the input rather than the number of operations in a forward pass. For the variable projected super-resolution problem, the memory consumption due to forward mode should remain constant regardless the number of CG iterations. In our experiments we observed that memory usage increases with the number of optimization steps which can be attributed to an implementation issue in PyTorch. 

The computational overhead of using CG to compute the projection still poses a difficulty and its of great interest to reduce the number of iterations while still obtaining quality reconstructions. We explore two options to reduce computational overhead: differentiating through only a portion of CG iterations, but solving the system near exactly, and differentiating through every iteration of CG, but not iterating CG to convergence when computing the projection.

\subsubsection{Differentiation of a Portion of CG Iterations}

\begin{figure}[]
    \centering
    \includegraphics[width=.6\linewidth]{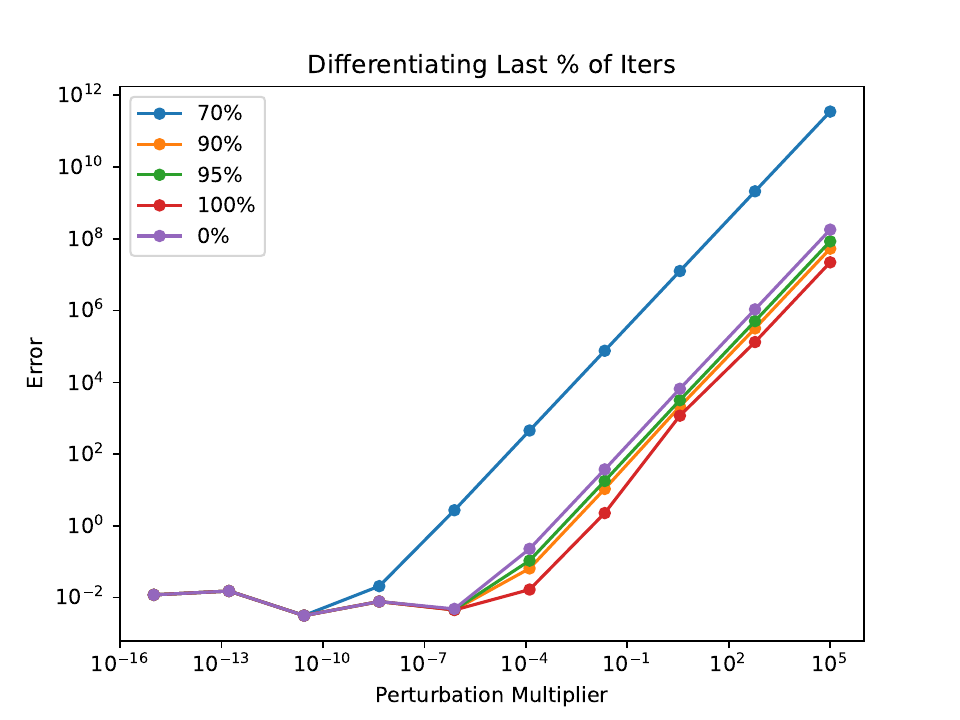}
    \caption{Jacobian convergence plot for the residual function $\bfr ( \bfw ) = \sfA_{\lambda_{\bff}} ( \bfw ) \bff^{(0)} ( \bfw ) - \bfb$. In this case, $\text{Error} = \| \bfr (\bfw+h\bfv) - \bfr (\bfw) - h \frac {\partial bfr (\bfw)}{\partial \bfw} \bfv \|,$ where $h$ is the perturbation multiplier, $v$ is the direction of the perturbation, and our Jacobian approximation $\frac {\partial bfr (\bfw)}{\partial \bfw}$ varies in accuracy depending on the percentage of final iterations of CG that are differentiated.}
    \label{fig:diff-last-jacobians}
\end{figure}

Differentiating through the later iterations of CG proves essential for accurate computation. In these experiments, a total of $200$ iterations of CG are used in each optimization step, but the percentage of the final iterations that are differentiated through vary. Figure \ref{fig:diff-last-jacobians} presents convergence plots for Jacobians obtained via AD, comparing various percentages of CG iterations through which differentiation is applied. Interestingly, differentiating through the final 70\% of iterations leads to poorer convergence than not differentiating through any iterations at all. However, when differentiating through the entire set of iterations (100\%), the convergence rate surpasses that of the 90\% and 95\% cases.

\begin{figure}[]
    \centering
    \includegraphics[width=\linewidth]{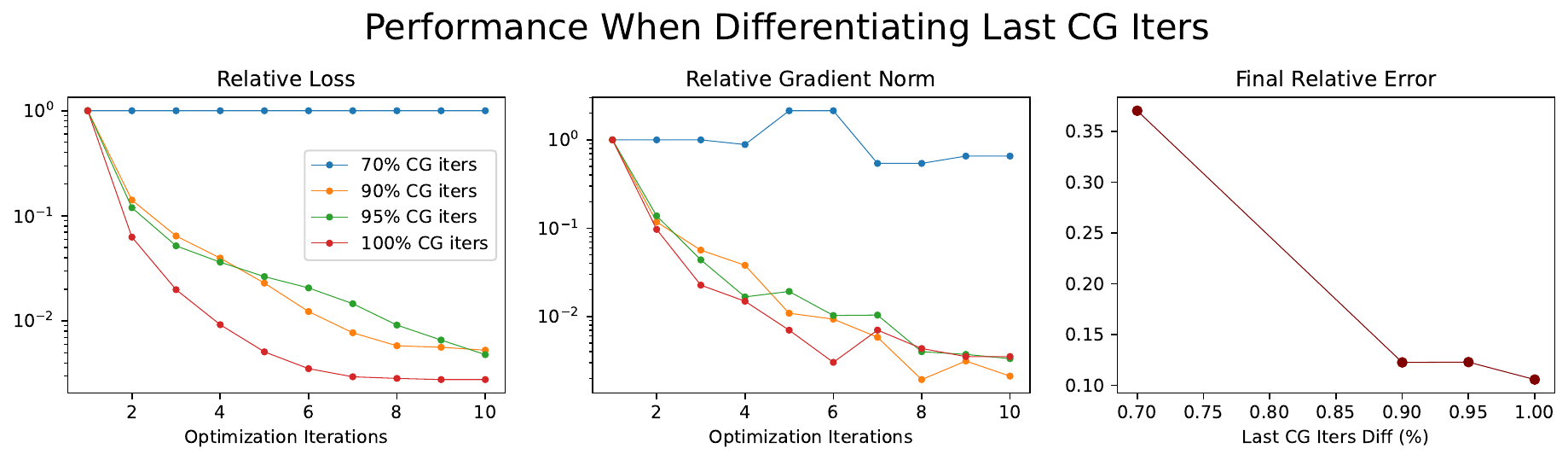}
    \caption{The performance of 200 iterations of variable-projected Gauss-Newton on the problem
from Figure \ref{fig:inexact_supperes_test} for various percentages of final CG iterations differentiated. The
final relative error is computed by solving the least squares problem exactly for the final reference
image prediction and measuring the difference from the true reference used to create the data.}
    \label{fig:diff-last-performance}
\end{figure}

\begin{figure}[]
    \centering
    \includegraphics[width=\linewidth]{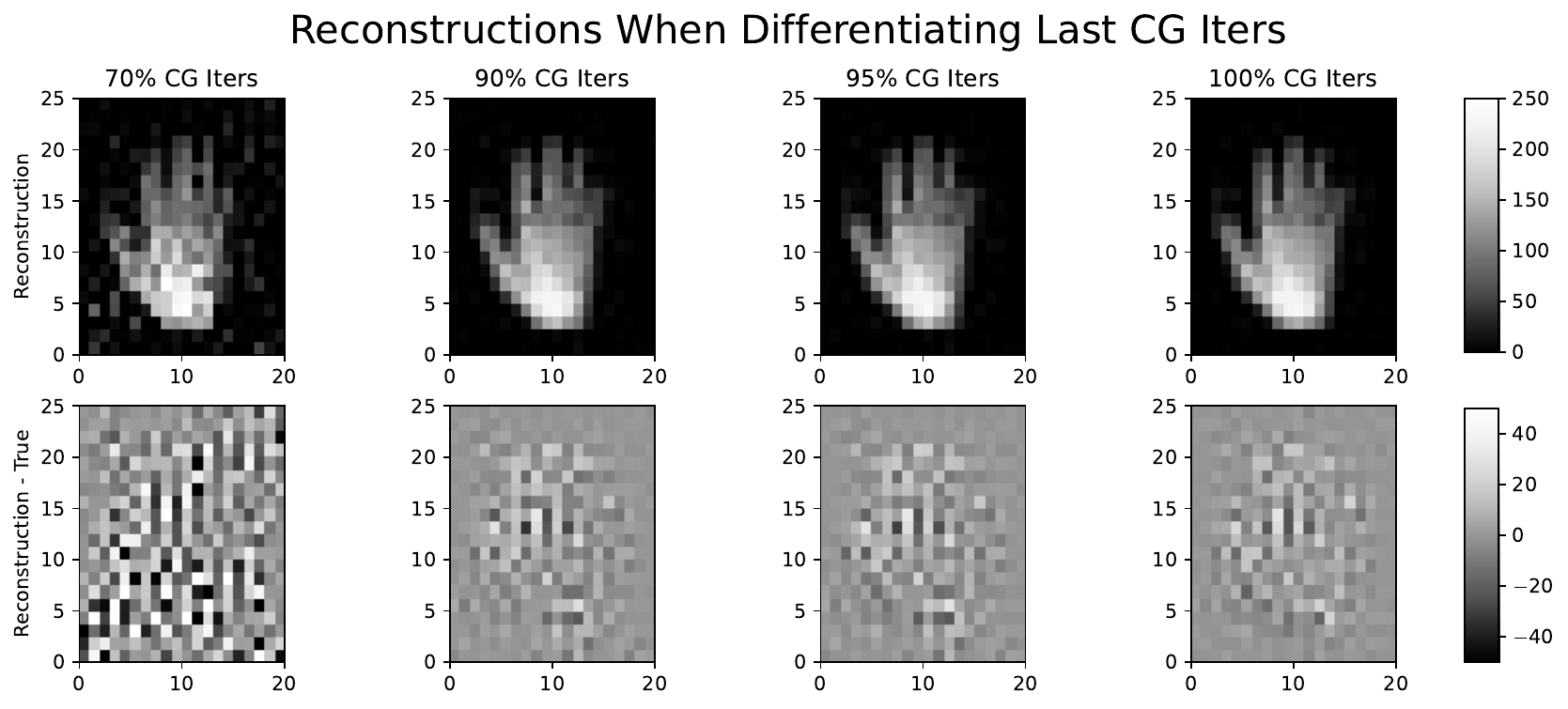}
    \caption{The final reference images and reconstruction errors of variable projected super resolution with inexact projections. The images are very similar, but some differences in the residuals
can be seen.}
    \label{fig:diff-last-reconstructions}
\end{figure}

Our results shown in Figure \ref{fig:diff-last-performance} demonstrate that differentiating through more than $70$\% is necessary. Differentiating through $100$\% of the iterations of CG is the best in terms of relative loss and gradient norm reduction. The final reconstruction qualities additionally indicate that $100$\% is the best, but the visual difference in error is small after $90$\% as seen in Figure \ref{fig:diff-last-reconstructions}. We also tested differentiating through the initial iterations of CG, but the performance was worse than using the later iterations. Therefore, these results are not included in this work. 

\subsubsection{Using Inexact Projections}

We test how reconstruction quality varies when using an inexact least squares solution for our projection. We consider the problem shown in Figure \ref{fig:inexact_supperes_test}, which we solve using a variable projected Gauss-Newton method. We automatically differentiate through every iteration of CG, but run a fixed number of CG iterations in each evaluation of the objective function.

\begin{figure}[]
    \centering
    \includegraphics[width=\linewidth]{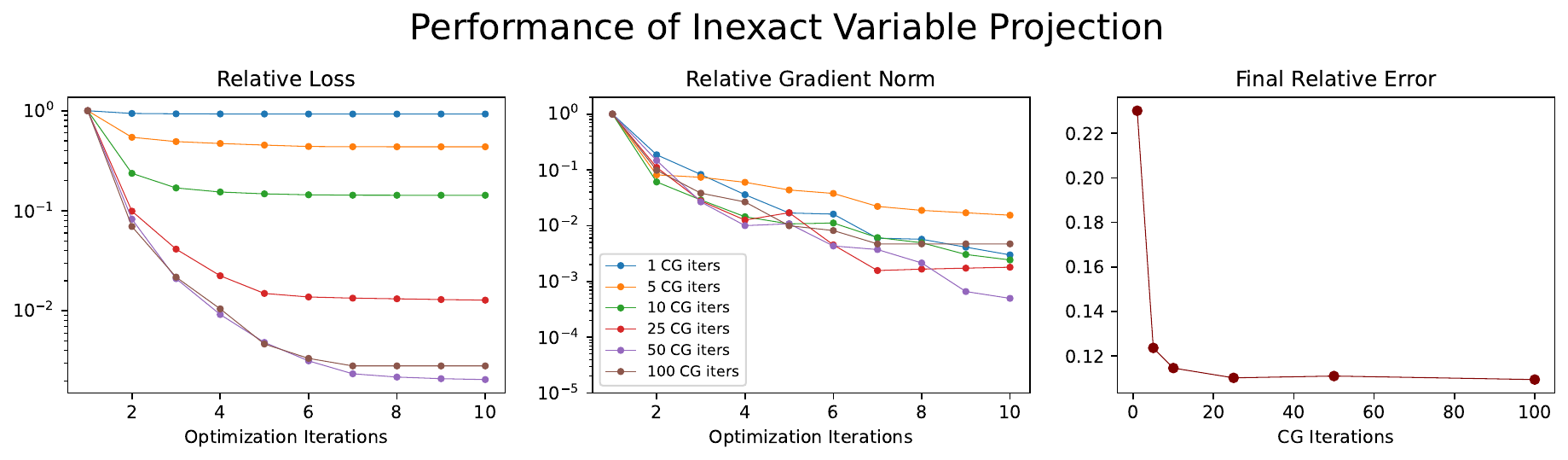}
    \caption{The performance of $10$ iterations of variable-projected Gauss-Newton on the problem from Figure \ref{fig:inexact_supperes_test} for various numbers of Conjugate Gradient Iterations used in each projection. The final relative error is computed by solving the least squares problem exactly for the final reference image prediction and measuring the difference from the true reference used to create the data.}
    \label{fig:inexact_supperes_results}
\end{figure}

\begin{figure}[]
    \centering
    \includegraphics[width=\linewidth]{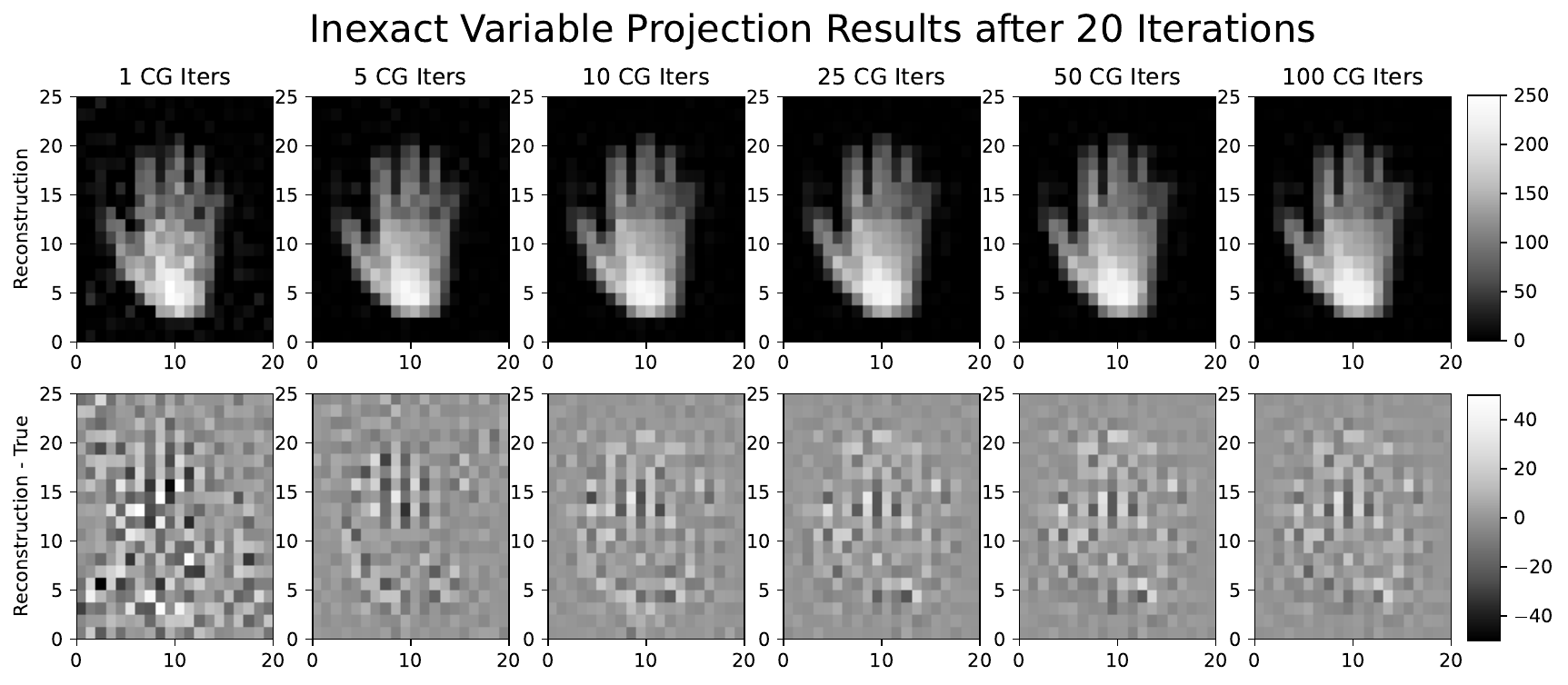}
    \caption{The final reference images and reconstruction errors of variable projected super resolution with inexact projections. The images are very similar, but some differences in the residuals can be seen.}
    \label{fig:inexact_supperes_final}
\end{figure}

Our results (Figure \ref{fig:inexact_supperes_results}) show that there are diminishing returns for using more exact projections in each iteration. Using $50$ or $100$ iterations performed the best in terms of relative loss and gradient norm reduction. The final recovered reference images (Figure \ref{fig:inexact_supperes_final}) show that the best reconstruction qualitatively is that using $100$ CG iterations, but $10$, $25$, and $50$ are close to the same reconstruction quality and the visual difference in reconstructions and error images is small.

\section{Conclusions}

In the affine registration problem, automatic differentiation has opened up the ability to directly compute Hessians, allowing for the application of Newton's method, which provides faster convergence than the previously used approaches. Furthermore, having access to exact Hessians allows for the implementation of a predictor-corrector method, which we've found to succeed in some problems where the standard multi-scale approach fails.

In the problem of super resolution, we showed that computing the full Jacobian of a variable projected objective function improves both the speed of convergence of the optimization and the final reconstruction quality from the previous approach of not differentiating through the projection, though it introduces computational overhead. Through experiments on a small super resolution problem, we found that tracking the computations of at least $70$\% of the final iterations of the CG method was necessary to see the benefits of including the Jacobian of the projection. Moreover, we found that using inexact projections during each iteration of optimization maintains high reconstruction quality with fewer CG iterations, providing a balance between computational efficiency and accuracy. These results suggest that leveraging inexact projections can mitigate computational costs while maintaining the advantage in results that the AD computed Jacobian of the projection provide.

Using AD tools, we have been able to simply compute derivatives, opening up methods in image registration that would be difficult to compute manually. The ability to automatically differentiate through complex processes, such as inexact iterative least squares solves, and registration objective functions, has proven its use in developing better super resolution solutions and second order optimization methods. These advancements demonstrate the potential of AD to further progress in image registration. Future work could explore scaling these methods to larger datasets and non-parametric transformations, as well as better optimizing memory usage when computing the derivative of inexact iterative least squares solves.

\section*{Acknowledgments}
This work was supported by the US National Science Foundation (NSF) award DMS-2349534. Any opinions, findings, and conclusions or recommendations expressed in this material are those of the authors and do not necessarily reflect the views of the NSF. We want to especially thank our mentor Dr. Lars Ruthotto and the other mentors of Emory's 2024 REU for Computational Mathematics and Data Science for making this possible.

\bibliographystyle{siamplain}
\bibliography{references}

\begin{thebibliography}{10}

\bibitem{chen2018neural}
{\sc R.~T. Chen, Y.~Rubanova, J.~Bettencourt, and D.~K. Duvenaud}, {\em Neural ordinary differential equations}, Advances in neural information processing systems, 31 (2018).

\bibitem{CGdiff}
{\sc B.~Christianson}, {\em Differentiating through conjugate gradient}, Optimization Methods and Software, 33 (2018), pp.~988--994, \url{https://doi.org/10.1080/10556788.2018.1425862}, \url{https://doi.org/10.1080/10556788.2018.1425862}, \url{https://arxiv.org/abs/https://doi.org/10.1080/10556788.2018.1425862}.

\bibitem{chung2006numerical}
{\sc J.~Chung, E.~Haber, and J.~Nagy}, {\em Numerical methods for coupled super-resolution}, Inverse Problems, 22 (2006), p.~1261, \url{https://doi.org/10.1088/0266-5611/22/4/009}, \url{https://dx.doi.org/10.1088/0266-5611/22/4/009}.

\bibitem{conn2000trust}
{\sc A.~R. Conn, N.~I. Gould, and P.~L. Toint}, {\em Trust region methods}, SIAM, 2000.

\bibitem{osti_876373}
{\sc D.~M. Dunlavy and D.~P. O'Leary}, {\em Homotopy optimization methods for global optimization.}, tech. report, Sandia National Laboratories (SNL), Albuquerque, NM, and Livermore, CA (United States), 12 2005, \url{https://doi.org/10.2172/876373}, \url{https://www.osti.gov/biblio/876373}.

\bibitem{español2024convergenceanalysisvariableprojection}
{\sc M.~I. Español and G.~Jeronimo}, {\em Convergence analysis of a variable projection method for regularized separable nonlinear inverse problems}, 2024, \url{https://arxiv.org/abs/2402.08568}, \url{https://arxiv.org/abs/2402.08568}.

\bibitem{Fischer_2008}
{\sc B.~Fischer and J.~Modersitzki}, {\em Ill-posed medicine—an introduction to image registration}, Inverse Problems, 24 (2008), p.~034008, \url{https://doi.org/10.1088/0266-5611/24/3/034008}, \url{https://dx.doi.org/10.1088/0266-5611/24/3/034008}.

\bibitem{golub2003separable}
{\sc G.~Golub and V.~Pereyra}, {\em Separable nonlinear least squares: the variable projection method and its applications}, Inverse problems, 19 (2003), p.~R1.

\bibitem{varpro}
{\sc G.~H. Golub and V.~Pereyra}, {\em The differentiation of pseudo-inverses and nonlinear least squares problems whose variables separate}, SIAM Journal on numerical analysis, 10 (1973), pp.~413--432.

\bibitem{CGdiff1}
{\sc S.~Gratton, D.~Titley-Peloquin, P.~Toint, and J.~T. Ilunga}, {\em Differentiating the method of conjugate gradients}, SIAM Journal on Matrix Analysis and Applications, 35 (2014), pp.~110--126, \url{https://doi.org/10.1137/120889848}, \url{https://doi.org/10.1137/120889848}, \url{https://arxiv.org/abs/https://doi.org/10.1137/120889848}.

\bibitem{Herring2018}
{\sc J.~L. Herring, J.~G. Nagy, and L.~Ruthotto}, {\em Lap: A linearize and project method for solving inverse problems with coupled variables}, Sampling Theory in Signal and Image Processing, 17 (2018), p.~127–151, \url{https://doi.org/10.1007/bf03549661}, \url{http://dx.doi.org/10.1007/BF03549661}.

\bibitem{functorch2021}
{\sc R.~Z. Horace~He}, {\em functorch: Jax-like composable function transforms for pytorch}.
\newblock \url{https://github.com/pytorch/functorch}, 2021.

\bibitem{otherdiff}
{\sc P.~Hovland and J.~Hückelheim}, {\em Differentiating through linear solvers}, 2024, \url{https://arxiv.org/abs/2404.17039}, \url{https://arxiv.org/abs/2404.17039}.

\bibitem{margossian2019review}
{\sc C.~C. Margossian}, {\em A review of automatic differentiation and its efficient implementation}, Wiley interdisciplinary reviews: data mining and knowledge discovery, 9 (2019), p.~e1305.

\bibitem{Modersitzki2009}
{\sc J.~Modersitzki}, {\em {FAIR: flexible algorithms for image registration}}, vol.~6 of Society for Industrial and Applied Mathematics (SIAM), Philadelphia, PA, Society for Industrial and Applied Mathematics (SIAM), Philadelphia, PA, 2009, \url{https://doi.org/10.1137/1.9780898718843}.

\bibitem{nnvarpro}
{\sc E.~Newman, L.~Ruthotto, J.~Hart, and B.~van Bloemen~Waanders}, {\em Train like a (var)pro: Efficient training of neural networks with variable projection}, SIAM Journal on Mathematics of Data Science, 3 (2021), pp.~1041--1066, \url{https://doi.org/10.1137/20M1359511}, \url{https://doi.org/10.1137/20M1359511}, \url{https://arxiv.org/abs/https://doi.org/10.1137/20M1359511}.

\bibitem{NoceWrig06}
{\sc J.~Nocedal and S.~J. Wright}, {\em Numerical Optimization}, Springer, New York, NY, USA, 2e~ed., 2006.

\bibitem{PyTorch}
{\sc A.~Paszke, S.~Gross, F.~Massa, A.~Lerer, J.~Bradbury, G.~Chanan, T.~Killeen, Z.~Lin, N.~Gimelshein, L.~Antiga, A.~Desmaison, A.~Kopf, E.~Yang, Z.~DeVito, M.~Raison, A.~Tejani, S.~Chilamkurthy, B.~Steiner, L.~Fang, J.~Bai, and S.~Chintala}, {\em Pytorch: An imperative style, high-performance deep learning library}, in Advances in Neural Information Processing Systems 32, Curran Associates, Inc., 2019, pp.~8024--8035, \url{http://papers.neurips.cc/paper/9015-pytorch-an-imperative-style-high-performance-deep-learning-library.pdf}.

\bibitem{schoenholz2020jax}
{\sc S.~S. Schoenholz and E.~D. Cubuk}, {\em {\{}JAX{\}} {\{}md{\}}: End-to-end differentiable, hardware accelerated, molecular dynamics in pure python}, 2020, \url{https://openreview.net/forum?id=r1xMnCNYvB}.

\bibitem{SUN2024103249}
{\sc S.~Sun, K.~Han, C.~You, H.~Tang, D.~Kong, J.~Naushad, X.~Yan, H.~Ma, P.~Khosravi, J.~S. Duncan, and X.~Xie}, {\em Medical image registration via neural fields}, Medical Image Analysis, 97 (2024), p.~103249, \url{https://doi.org/https://doi.org/10.1016/j.media.2024.103249}, \url{https://www.sciencedirect.com/science/article/pii/S1361841524001749}.

\bibitem{2020SciPy-NMeth}
{\sc P.~Virtanen, R.~Gommers, T.~E. Oliphant, M.~Haberland, T.~Reddy, D.~Cournapeau, E.~Burovski, P.~Peterson, W.~Weckesser, J.~Bright, S.~J. {van der Walt}, M.~Brett, J.~Wilson, K.~J. Millman, N.~Mayorov, A.~R.~J. Nelson, E.~Jones, R.~Kern, E.~Larson, C.~J. Carey, {\.I}.~Polat, Y.~Feng, E.~W. Moore, J.~{VanderPlas}, D.~Laxalde, J.~Perktold, R.~Cimrman, I.~Henriksen, E.~A. Quintero, C.~R. Harris, A.~M. Archibald, A.~H. Ribeiro, F.~Pedregosa, P.~{van Mulbregt}, and {SciPy 1.0 Contributors}}, {\em {{SciPy} 1.0: Fundamental Algorithms for Scientific Computing in Python}}, Nature Methods, 17 (2020), pp.~261--272, \url{https://doi.org/10.1038/s41592-019-0686-2}.

\bibitem{repo}
{\sc W.~Watson, C.~Cherry, R.~Lang, and L.~Ruthotto}, {\em Imgregpytorchproject}.
\newblock \url{https://github.com/wdwatson2/ImgRegPytorchProject}, 2024.
\newblock GitHub repository.

\bibitem{wu2022nodeo}
{\sc Y.~Wu, T.~Z. Jiahao, J.~Wang, P.~A. Yushkevich, M.~A. Hsieh, and J.~C. Gee}, {\em Nodeo: A neural ordinary differential equation based optimization framework for deformable image registration}, in Proceedings of the IEEE/CVF conference on computer vision and pattern recognition, 2022, pp.~20804--20813.

\end{thebibliography}

\end{document}